\documentclass{article}
\usepackage{graphicx} 
\usepackage{hyperref}
\usepackage{amssymb}
\usepackage{amsmath}

\newtheorem{theorem}{Theorem}

\newtheorem{corollary}[theorem]{Corollary}

\title{The Maximum of the Volume of a Part of a Cevian Simplex}
\author{Zamina Guliyeva, Yagub Aliyev}
\date{January 2026}

\begin{document}

\maketitle

\begin{abstract}
The cevians passing through a point in a simplex create a cevian simplex, which is divided by these cevians into smaller simplices. We consider the problem about the maximum of the ratio of the sum of the volumes of some of these smaller simplices by the volume of the reference simplex. The special case of tetrahedron is given as an example.

\end{abstract}

Given an $n$-simplex with vertices $A_1$, $A_2,\ldots$, $A_{n+1}$ ($n\ge2$) and any point $M$ in its interior, let $N_i$ be the points on the hyperplane opposite $A_i$ such that $A_i$, $M$, and $N_i$ are collinear. Let $V=[A_1A_2\ldots A_{n+1}]$ and $$V_i=[MN_1N_2\ldots N_{i-1} N_{i+1}\ldots N_{n+1}]\ (1\le i\le n+1),$$where the square brackets indicate the volume of the $n$-simplex between them. We want to find the point $M$ for which $\frac{V_1+V_2+\cdots +V_k}{V}$ ($1\le k\le n+1$) is maximal and find this maximum. The case $k=n+1$ was considered in \cite{reut}, \cite{klam}, \cite{mit}, p. 502. The cases $k=1$ and $k=n$ were considered in \cite{ali}. 
In the current paper we will consider the remaining cases $1<k<n$.

Let $\overline\lambda=(\lambda_1,\lambda_2,\ldots,\lambda_{n+1})$ be the barycentric coordinates of $M$ with respect to the simplex $A_1A_2\ldots A_{n+1}$. Then $\lambda_1+{\lambda_2}+\cdots+\lambda_{n+1}=1$. Suppose that $R_i=|MA_i|$, and $s_i=|MN_i|$. Since $\lambda_i=\frac{s_i}{R_i+s_i}$, $\frac{s_i}{R_i}=\frac{\lambda_i}{1-\lambda_i}$ (see \cite{moeb}; \cite{samet}; \cite{mark}; \cite{balk}, p. 124-126). Therefore, 
$$V_i=[MA_1A_2\ldots A_{i-1}A_{i+1}\ldots A_{n+1}]\cdot\prod_{j\ne i}\frac{s_j}{R_j}=V\cdot (1-\lambda_{i})\cdot \prod_{j=1}^{n+1}\frac{\lambda_j}{1-\lambda_j}.$$
The same result, up to a sign, can be obtained by calculating $V_i$ using the coordinates of the vertices of $MN_1N_2\ldots N_{i-1} N_{i+1}\ldots N_{n+1}$
$$
\frac{V_i}{V} =\pm \begin{vmatrix}
    \lambda_1 & \lambda_2 & \cdots & \lambda_{n+1} \\
    0 & \frac{\lambda_2}{1-\lambda_1} & \cdots & \frac{\lambda_{n+1}}{1-\lambda_1} \\
    \frac{\lambda_1}{1-\lambda_2} & 0 & \cdots & \frac{\lambda_{n+1}}{1-\lambda_2} \\
    \cdots & \cdots & \cdots & \cdots\\
    \frac{\lambda_1}{1-\lambda_{n+1}} &\frac{\lambda_2}{1-\lambda_{n+1}} & \cdots &0\\
  \end{vmatrix},
$$
where the row corresponding to the point $N_i$ with barycentric coordinates
$$
\left(\frac{\lambda_1}{1-\lambda_i},\ldots,\frac{\lambda_{i-1}}{1-\lambda_i},0,\frac{\lambda_{i+1}}{1-\lambda_i},\ldots,\frac{\lambda_{n+1}}{1-\lambda_i}  \right),
$$
is omitted.
Consequently,
$$
V_1+V_2+\cdots +V_k=V\cdot(k-\lambda_1-\lambda_2-\cdots-\lambda_k)\cdot \prod_{j=1}^{n+1}\frac{\lambda_j}{1-\lambda_j}.
$$
In particular, if $k=n+1$, then we obtain $
V_1+V_2+\cdots +V_{n+1}=V\cdot F(\overline\lambda),
$ where 
$$
F(\lambda_1,\lambda_2,\ldots,\lambda_{n+1})= n\cdot \prod_{j=1}^{n+1}\frac{\lambda_j}{1-\lambda_j}.
$$
The function $F_k(\lambda_1,\lambda_2,\ldots,\lambda_{n+1})$ is continuous whenever $\lambda_i\ne 1$ for all $i=1,\ldots,n+1$. In order to make it continuous at the vertices $A_i$ ($i=1,\ldots,n+1$) of the simplex, where $\lambda_i=1$, we need to define $F(\lambda_1,\lambda_2,\ldots,\lambda_{n+1})=0$. After redefining the function $F(\lambda_1,\lambda_2,\ldots,\lambda_{n+1})$, the continuity at the vertex $A_i$ follows from the fact that if $\lambda_i\rightarrow 1$, then $\lambda_j\rightarrow 0$ for all $j\ne i$, and therefore $\prod_{j\ne i}\lambda_j=o\left(1-\lambda_i\right)$. Indeed,
$$
\prod_{j\ne i}\lambda_j=\frac{1}{n}\cdot \sum_{k=1}^{n}\prod_{j\ne i}\lambda_j=\frac{1}{n}\cdot \sum_{j\ne i}o(\lambda_j)=o\left(\sum_{j\ne i}\lambda_j\right)=o\left(1-\lambda_i\right).
$$
For $1<k<n$, we can write $
\frac{V_1+V_2+\cdots +V_k}{V}= F_k(\lambda_1,\lambda_2,\ldots,\lambda_{n+1}),
$ where $$F_k(\lambda_1,\lambda_2,\ldots,\lambda_{n+1})=G_k(\lambda_1,\lambda_2,\ldots,\lambda_{k})\cdot H_k(\lambda_{k+1},\lambda_{k+2},\ldots,\lambda_{n+1}),$$
$$
G_k(\lambda_1,\lambda_2,\ldots,\lambda_{k})= (k-\lambda_1-\lambda_2-\cdots-\lambda_k)\cdot \prod_{j=1}^{k}\frac{\lambda_j}{1-\lambda_j},
$$
$$
H_k(\lambda_{k+1},\lambda_{k+2},\ldots,\lambda_{n+1})=\prod_{j=k+1}^{n+1}\frac{\lambda_j}{1-\lambda_j}.
$$
Since $F_k(\lambda_1,\lambda_2,\ldots,\lambda_{n+1})\le F(\lambda_1,\lambda_2,\ldots,\lambda_{n+1})$, we can redefine $F_k(\overline\lambda)=0$ at the vertices $A_i$ ($i=1,\ldots,n+1$) of the simplex to make it continuous at these vertices, too.
Now the function $F_k(\lambda_1,\lambda_2,\ldots,\lambda_{n+1})$ is continuous in the simplex $A_1A_2\ldots A_{n+1}$, including its boundary, and therefore by the classical result in analysis (see e.g. \cite{rud}, Theorem 4.16) $F_k(\overline\lambda)$ attains both maximum and minimum values in this simplex, which is a compact set.
At the boundary (vertices, edges, faces, etc.) $F_k(\overline\lambda)=0$, which
is the minimum. So, the maximum occurs in the interior of the simplex.
For fixed ${\lambda_i}+{\lambda_j}$ the product $\frac{\lambda_i}{1-\lambda_i}\cdot \frac{\lambda_j}{1-\lambda_j}$ reaches the maximum when $\lambda_i={\lambda_j}$. Therefore, for fixed $\sum_{i=1}^k{\lambda_i}=kx$ the function $G_k$ reaches the maximum when $\lambda_1={\lambda_2}=\cdots=\lambda_k=x$. Similarly, for fixed $\sum_{i=k+1}^{n+1}{\lambda_i}=my$, where $m=n+1-k>0$, the function $H_k$ reaches the maximum when $\lambda_{k+1}=\lambda_{k+2}=\cdots=\lambda_{n+1}=y$. Then $y=\frac{1-kx}{m}$, $1-y=\frac{m-1+kx}{m}$, and $$F_k(\lambda_1,\lambda_2,\ldots,\lambda_{n+1})=k\cdot(1-x)\cdot\left(\frac{x}{1-x}\right)^{k}\cdot \left(\frac{1-kx}{m-1+kx}\right)^m.$$
By studying the function $f(x)=k\cdot\frac{x^k}{\left(1-x\right)^{k-1}}\cdot \left(\frac{1-kx}{m-1+kx}\right)^m$ for the maximum in interval $\left[0,\frac{1}{k}\right]$ we find that
$$f'(x)=k\cdot\frac{x^{k-1}}{\left(1-x\right)^{k}}\cdot \frac{(1-kx)^{m-1}}{(m-1+kx)^{m+1}}\cdot q(x),$$
where
$$
q(x)=\left(k-x\right)\left(1-kx\right) \left(m-1+k x\right)-km^{2}x\left(1-x\right).
$$
Note that $q(-\infty)=-\infty$, $q(0)=k(m-1)>0$, $q\left(\frac{1}{k}\right)=-\frac{m^2(k-1)}{k}<0$, and $q(+\infty)=+\infty$. Therefore, the cubic $q(x)$ has exactly one zero in each of the intervals $(-\infty,0)$, $\left(0,\frac{1}{k}\right)$, and $\left(\frac{1}{k},+\infty \right)$.
Let us denote by $\theta$ the zero of $q(x)$ in $\left(0,\frac{1}{k}\right)$. Consequently, $$f(x)\le f \left(\theta\right)=k\cdot\frac{\theta^k}{\left(1-\theta\right)^{k-1}}\cdot \left(\frac{1-k\theta}{m-1+k\theta}\right)^{n+1-k}.$$
Thus we proved the following result.
\begin{theorem}
    Let $A_1A_2\ldots A_{n+1}$ $(n\ge2)$ be an $n$-simplex and $M$ be any point in its interior. Let $N_i$ be the points on the hyperplane opposite $A_i$ such that $A_i$, $M$, and $N_i$ are collinear. Let $V$ be the volume of $A_1A_2\ldots A_{n+1}$ and let $V_i$ be the volume of $MN_1N_2\ldots N_{i-1} N_{i+1}\ldots N_{n+1}$ $(1\le i\le n+1)$. Then for $1<k<n$
$$
\frac{V_1+V_2+\cdots +V_k}{V}\le \frac{k\theta^k}{\left(1-\theta\right)^{k-1}}\cdot \left(\frac{1-k\theta}{m-1+k\theta}\right)^{m},
$$ 
where $m=n+1-k$, and $\theta$ is the zero of the cubic $$q(x)=k^{2}x^{3}
- k(k^{2}-m^{2}-m+2)x^{2}
+ (2k^2-k^2m-km^{2}-m+1)x
+ k(m-1)
$$
in $\left(0,\frac{1}{k}\right)$. Equality happens at the point with barycentric coordinates $\lambda_1={\lambda_2}=\cdots=\lambda_k=\theta$,  $\lambda_{k+1}=\lambda_{k+2}=\cdots=\lambda_{n+1}=\frac{1-k\theta}{m}$.
\end{theorem}
\textbf{Remark 1.} If $n=3$, $k=2$, and $m=2$, then we obtain the following result in tetrahedron geometry (see Fig. 1). The calculations are done in Maple 2024.
\begin{corollary}
Let $A_1A_2A_3A_{4}$ be a tetrahedron and $M$ be any point in its interior. Let $N_i$ be the points on the plane opposite $A_i$ such that $A_i$, $M$, and $N_i$ are collinear. Let $V$, $V_1$, and $V_2$ be the volumes of tetrahedra $A_1A_2A_3A_{4}$, $MN_2N_3N_{4}$, and $MN_1N_3N_{4}$, respectively. Then
$$
\frac{V_1+V_2}{V}\le \frac{2\theta_0^{2}(1-2\theta_0)}{3} \approx 0.018800151866697,
$$ 
where $\theta_0=\sqrt{3}\, \cos\! \left(\frac{\arctan\! \left(\frac{\sqrt{23}}{2}\right)}{3}+\frac{\pi}{3}\right)\approx 0.22745204561175$ is the zero of the cubic $q(x)=4x^{3}-9x+2$
in $\left(0,\frac{1}{2}\right)$.
\end{corollary}
\begin{figure}[htb]
\begin{center}
\includegraphics{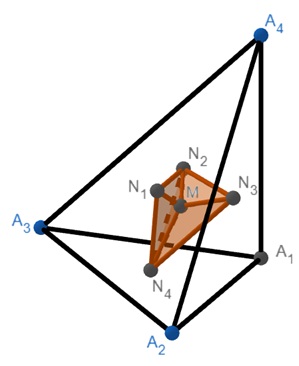}
\caption{\url{https://www.geogebra.org/3d/htesqury}}
\end{center}
\end{figure}
\textbf{Remark 2.} If $k=0$, $k=1$, $m=0$, or $m=1$, then the cubic $q(x)$ has a rational zero. It would be interesting to find $k>1$ and $m>1$ such that the cubic $q(x)$ has a rational zero or prove that such $k$ and $m$ do not exist.


\begin{thebibliography}{HD82}
\bibitem{ali} Yagub N. Aliyev, The maximum of the volume of a cevian simplex and its parts, The American Mathematical Monthly, accepted.
\end{thebibliography}
\end{document}